\documentclass[12pt]{amsart}
\usepackage[dvips]{graphics}
\usepackage{amssymb, amsmath}

\newtheorem{theorem}{Theorem}[section]
\newtheorem{proposition}[theorem]{Proposition}

\newtheorem{lemma}[theorem]{Lemma}

\theoremstyle{definition}

\theoremstyle{remark}
\newtheorem{remark}[theorem]{Remark}
\newtheorem{remarks}[theorem]{Remarks}

\def\varle{\leqslant}

\newcommand{\CM}{{\mathcal M}}

\newcommand{\CO}{{\mathcal O}}
\newcommand{\CP}{{\mathcal P}}

\newcommand{\CS}{{\mathcal S}}

\newcommand{\CW}{{\mathcal W}}

\newcommand{\FM}{{\bf M}}

\newcommand{\FP}{{\bf P}}

\newcommand{\fh}{{{\mathfrak h}}} 
 
\newcommand{\fp}{{{\mathfrak p}}} 
\newcommand{\fg}{{{\mathfrak g}}} 
\newcommand{\fb}{{{\mathfrak b}}} 
 
\newcommand{\fn}{{{\mathfrak n}}} 
\newcommand{\fm}{{{\mathfrak m}}}

\newcommand{\fhd}{\fh^\star}

\newcommand{\DC}{{\mathbb C}}

\newcommand{\DZ}{{\mathbb Z}}
\newcommand{\DK}{{\mathbb K}}

\newcommand{\DN}{{\mathbb N}}

\newcommand{\DV}{{\mathbb V}}

\newcommand{\ch}{{\operatorname{ch}}}

\newcommand{\End}{{\operatorname{End}}}

\newcommand{\Hom}{{\operatorname{Hom}}}
\newcommand{\id}{{\operatorname{id}}}

\newcommand{\Ind}{{\operatorname{Ind}}}
\newcommand{\Res}{{\operatorname{Res}}}
\newcommand{\catmod}{{\operatorname{-mod}}}
\newcommand{\catmodr}{{\operatorname{mod-}}}

\newcommand{\Stab}{{\operatorname{Stab}}}
\newcommand{\sur}{\mbox{$\rightarrow\!\!\!\!\!\rightarrow$}}

\begin{document}

\pagenumbering{arabic}
\title[The combinatorics of $\CO$]{The Combinatorics of Category $\CO$ over symmetrizable Kac-Moody Algebras}

\begin{abstract}
We show that the structure of blocks outside the critical hyperplanes of category $\CO$ over any symmetrizable Kac-Moody algebra depends only on the corresponding integral Weyl group and its action on the parameters of the Verma modules by giving a combinatorial description of the projective objects. As an application we derive the Kazhdan-Lusztig conjecture for non--integral blocks from the integral case in finite and affine situations.

\end{abstract}
\author[]{Peter Fiebig${}^\ast$}
\thanks{${}^\ast$ supported by the German Academic Exchange Service (DAAD)}

\maketitle

\section{Introduction}
One of the most prominent categories of representations of triangulated Lie algebras is the category $\CO$, originally introduced by Bernstein, Gelfand and Gelfand.  It provides the natural framework for the study of highest weight modules. In \cite{Soe1} Soergel showed that in the case of a finite dimensional, semisimple complex Lie algebra one can give a description of the categorical structure of $\CO$ in terms of the corresponding Weyl group. This article provides the corresponding statement for symmetrizable Kac-Moody algebras. 

Let us first review Soergel's result and the main ideas in its proof. Let $\fg$ be a finite dimensional semisimple complex Lie algebra and choose a Borel subalgebra $\fb\subset\fg$ and a Cartan subalgebra $\fh\subset \fb$. The simple objects $L(\lambda)$ of the corresponding category $\CO$ are then parametrized by elements $\lambda\in\fhd$. Let $\CO=\prod_\Lambda\CO_\Lambda$ be the decomposition of $\CO$ into indecomposable blocks.  We can identify a parameter $\Lambda$ with the set of simple objects in $\CO_\Lambda$, hence with a subset of $\fhd$. Then each of these subsets $\Lambda$ is an orbit of a subgroup $\CW(\Lambda)$ of the Weyl group $\CW$ of $\fg$ under the linear action shifted by the half sum of positive roots. In particular, $\Lambda$ is finite. There is a prefered choice of simple reflections $\CS(\Lambda)\subset\CW(\Lambda)$ and $(\CW(\Lambda),\CS(\Lambda))$ forms a Coxeter system. Soergel proved that  $\CO_\Lambda$ depends only on the isomorphism class of this system and on its action on $\Lambda$. More precisely, if $\fh\subset\fb\subset\fg$ and $\fh^\prime\subset\fb^\prime\subset\fg^\prime$ are two semisimple Lie algebras together with Cartan and Borel subalgebras and $\CO_\Lambda$ and $÷\CO^\prime_{\Lambda^\prime}$ are two blocks of the corresponding categories such that there is an isomorphism $(\CW(\Lambda),\CS(\Lambda))\cong(\CW^\prime(\Lambda^\prime),\CS^\prime(\Lambda^\prime))$ of the associated Coxeter systems which induces a bijection $\Lambda\cong\Lambda^\prime$ of orbits, then there is an equivalence $\CO_\Lambda\cong\CO^\prime_{\Lambda^\prime}$ of categories. 

The proof consists of a combinatorial description of the structure of a block $\CO_\Lambda$ in terms of  $(\CW(\Lambda), \CS(\Lambda))$ and its action on $\Lambda$.
The main ideas are the following. First, let $\lambda\in\Lambda$ be the dominant weight. For any $w\in\CW(\Lambda)$ let $P(w.\lambda)$ be the projective cover of  $L(w.\lambda)$. Then $P=\bigoplus P(\bar w.\lambda)$, where the sum is taken over $\bar w\in\CW(\Lambda)/\Stab(\lambda)$, is a faithful small projective, i.e.ø
$\Hom(P,\cdot)\colon\CO_\Lambda\to\catmodr\End(P)$ is an equivalence of categories. Soergel gave a combinatorial, though not explicit description of $\End(P)$.  Let us assume that $\Lambda$ is regular, i.e.ø $\CW(\Lambda)$ acts faithfully on $\Lambda$. Let $w_0\in\CW(\Lambda)$ be the longest element, hence $P(w_0.\lambda)$ is the antidominant projective module. To the Coxeter system  $(\CW(\Lambda),\CS(\Lambda))$ one associates the commutative algebra of coinvariants $C=C(\CW(\Lambda),\CS(\Lambda))$ which comes equipped with an action of $\CW(\Lambda)$. For any $s\in\CS(\Lambda)$ let $C^s\subset C$ be the subalgebra of $s$--invariant elements. There is an isomorphism $C\cong\End(P(w_0.\lambda))$, hence we get a functor
$\DV:=\Hom(P(w_0.\lambda),\cdot)\colon \CO_{\Lambda}\to C\catmod$. The next two properties of $\DV$ are crucial.
\begin{list}{}{}
\item[(A)]  If $P$ and $P^\prime$ are projective in $\CO_\Lambda$, then $\DV$ induces an isomorphism $\Hom_{\CO}(P,P^\prime)\cong\Hom_{C}(\DV P,\DV P^\prime)$. 
\item[(B)] If $s\in\CS(\Lambda)$  and $\Theta_s\colon\CO_\Lambda\to\CO_\Lambda$ is the corresponding functor of ``translation through the $s$--wall", then $\DV\circ\Theta_s\cong C\otimes_{C^s}\DV$.
\end{list}
By (A) we have to describe $\DV P(w.\lambda)$ as a $C$--module. We do this by induction on the length of $w$ in the Coxeter system $(\CW(\Lambda), \CS(\Lambda))$. To start with, $\DV P(\lambda)\cong\DC$ is the unique simple quotient of $C$. Let  $w=s_1\cdots s_n$ be a reduced expression and choose any direct decomposition of $\Theta_{s_n}\cdots\Theta_{s_1} M(\lambda)$. Then $P(w.\lambda)$ is isomorphic to the indecomposable direct summand   that is not isomorphic to $P(w^\prime.\lambda)$ with $l(w^\prime)<l(w)$. By (A) and (B), $\DV P(w.\lambda)$ is the indecomposable direct summand of $C\otimes_{C^{s_n}}C\cdots\otimes_{C^{s_1}}\DC$ that is not isomorphic to $\DV P(w^\prime.\lambda)$ with $l(w^\prime)<l(w)$. So we inductively established the description of all $\DV P(w.\lambda)$ as modules over $C$ for regular $\Lambda$. If $\Lambda$ is not regular, one can give a desription of all projective objects using the regular case and a functor of ``translation on the walls''.

Now let $\fg$ be a symmetrizable Kac-Moody algebra. Then any parameter $\Lambda$ for the block decomposition of $\CO$ is again given by a subset of $\fhd$ and, if $\Lambda$ does not intersect the critical hyperplanes (i.e.ø the hyperplanes defined by an integrality condition on imaginary roots), then $\Lambda$ is again an orbit under a subgroup $\CW(\Lambda)$ of $\CW$, which is part of a Coxeter datum $(\CW(\Lambda),\CS(\Lambda))$. $\Lambda$ is not necessarily  finite.
 
From now on let $\Lambda$ be outside the critical hyperplanes. We prove a generalization of Soergel's result for the symmetrizable Kac-Moody case (Theorem \ref{equivalences of categories}). We follow the ideas explained above, though most statements need different proofs.

We distinguish two cases.
In the first case we assume that $\Lambda$ contains a dominant weight, i.e.ø a highest  weight under the usual partial order on $\fhd$. In the second case it contains an antidominant, i.e.ø lowest, weight. Suppose $\Lambda$ contains a dominant weight $\lambda\in\Lambda$. Then there exist projective covers $P(w.\lambda)$ of $L(w.\lambda)$ for any $w\in\CW(\Lambda)$ and the set $\CP_\Lambda:=\{P(\bar w.\lambda), \bar w\in\CW(\Lambda)/\Stab(\lambda)\}$   is a faithful set of small projectives in the sense of \cite{Mit}, i.e.ø the functor
\begin{eqnarray*}
\CO_\Lambda & \to & (\DC\catmod)^{\CP_\Lambda^{opp}} \\
M & \mapsto & \Hom(\cdot, M) 
\end{eqnarray*}
is an equivalence of categories, where $(\DC\catmod)^{\CP_\Lambda^{opp}}$ is the category of all additive functors $\CP_\Lambda^{opp}\to\DC\catmod$. Hence we have to describe $\CP_\Lambda$.  

If $\Lambda$ is infinite and contains an antidominant weight, there are no projective objects at all in $\CO_\Lambda$. However, there is a tilting equivalence $t\colon \CM\stackrel{\sim}\to\CM^{opp}$,
where $\CM\subset \CO$ is the full subcategory of modules which admit a Verma flag. It induces an equivalence of blocks $\CM_\Lambda\stackrel{\sim}\to\CM_{t(\Lambda)}$, where $t(\Lambda)$ will contain a dominant weight. Moreover, the structure of  $\CM(\Lambda)$ determines the structure of $\CO_\Lambda$, hence we have reduced the second case to the first.

We need to define the functor $\DV$. To start with, let $\lambda\in\Lambda$ be an antidominant weight. Then we construct the ``antidominant projective cover" $P^\infty(\lambda)$ as a certain limit of antidominant projective covers in truncated subcategories. Then we replace the algebra of coinvariants by the specialization $Z_\Lambda$ of the categorical (or Bernstein) center of a deformed version of $\CO_\Lambda$, which was calculated in \cite{Fie}. It only depends on $\CW(\Lambda)$ and its action on $\Lambda$. We show that there is a natural isomorphism $Z_\Lambda\cong\End(P^\infty(\lambda))$, so we get a functor
$\DV=\Hom(P^\infty(\lambda),\cdot)\colon \CO_\Lambda\to Z_\Lambda\catmod$. Restricting to the category of modules with Verma flag and using the tilting equivalence we analogously get a functor $\DV\colon\CM_\Lambda\to Z_\Lambda\catmod$ in the case when $\Lambda$ contains a dominant weight. Since every projective object in $\CO_\Lambda$ admits a Verma flag, this serves our purpose.

In \cite{Fie}, translation functors on $\CM_\Lambda$ were constructed and it was shown that they behave as in the finite dimensional situation, in particular they can be used to construct projective objects. So we have to verify statements (A) and (B). The original proofs in \cite{Soe1} do not work in the infinite setup, since $\DV$ is not a quotient functor. This might be remedied by certain limit and completion procedures, though we choose another way. We first prove both statements in the generic and subgeneric cases (Theorem \ref{struktursatz - subgeneric case} \& Theorem \ref{comb of translation - subgeneric}), i.e.ø when $\CW(\Lambda)$ is either trivial or contains only one reflection. This is done more or less by explicit verification. We then use the deformation construction as described in \cite{Fie} to reduce the general case to the generic or subgeneric cases (Theorem \ref{comb of translation} \& Theorem \ref{struktursatz}).

As an application we prove the Kazhdan-Lusztig conjecture for symmetrizable Kac-Moody algebras under the assumption that the associated Coxeter system  $(\CW(\Lambda),\CS(\Lambda))$ is of finite or affine type  by reduction to the integral case (Theorem \ref{Kazhdan-Lusztig conjecture}). As a byproduct of  the theory of translation functors we derive the uniqueness of embeddings of Verma modules outside the critical hyperplanes (Theorem \ref{Verma embeddings}), originally proved in the dominant case by Kashiwara-Tanisaki \cite{KT4}. 

I wish to thank Bernhard Keller for help with the categorical concepts used in this article and Wolfgang Soergel for many discussions. A major part of this work was done during a stay at the Mathematical Sciences Research Institute at Berkeley which was financed by the German Academic Exchange Service (DAAD). I am grateful  towards both institutions.

\section{Preliminaries}
Let $\fg$ be a complex symmetrizable Kac-Moody algebra, $\fb\subset\fg$ its Borel subalgebra and $\fh\subset\fb$ its Cartan subalgebra. Let $U=U(\fg)$, $B=U(\fb)$ and $S=U(\fh)=S(\fh)=\DC[\fhd]$ be the universal enveloping algebras. Let $\Pi\subset\Delta_+\subset\Delta\subset\fhd$ be the set of simple roots, the set of positive roots  and the set of roots of $\fg$ with respect to $\fh$ and $\fb$. Let $\Delta^{re}$ and $\Delta^{im}$ be the sets of real and imaginary roots. Define the usual partial order on $\fhd$ by setting $\lambda-\mu\ge0$ if and only if $\lambda-\mu\in\DN\,\Pi$. 

In this section we will quote the results and describe the methods that we will use in the following. Unless stated otherwise, the proofs can be found in \cite{Fie}. 

\subsection{Local deformation algebras}
We will call a commutative, associative, noetherian, unital $S$--algebra  which is a local domain a {\em local deformation algebra}. For any such algebra $T$ the $S$--structure is given by the structure morphism $\tau\colon S\to T$. We will be particularly interested in the following examples.

Let $R=S_{(\fh)}$ be the localization of $S$ at the maximal ideal generated by $\fh\subset S$, i.e.ø the localization at $0\in\fhd$. For any prime ideal $\fp\subset R$ let $R_\fp$ be the localization of $R$ at $\fp$ and let $\DK_\fp=R_\fp/R_\fp\fp$ be the corresponding residue field. The rings $R, R_\fp$ and $\DK_\fp$ are local deformation algebras. Note that as a special case the residue field $\DC$ of $R$ inherits an $S$--algebra structure, where the structure morphism $S\to \DC$ is given by evaluation at $0\in\fhd$.

A symmetrizable Kac-Moody algebra admits a non--degenerate, symmetric, invariant bilinear form $(\cdot,\cdot)\colon\fg\times\fg\to \DC$. It restricts to a non--degenerate form on $\fh\times\fh$ and induces a non--degenerate form on $\fhd\times\fhd$. Let $T$ be any local deformation algebra and let $\fh_T^\star:=\fhd\otimes_\DC T=\Hom_\DC(\fh,T)$. By T--bilinear extension we arrive at a non--degenerate form $(\cdot,\cdot)_T\colon\fh_T^\star\times\fh_T^\star\to T$. The restriction of the structure morphism $\tau\colon S\to T$ to $\fh\subset S$ is an element in $\fh_T^\star$. Moreover we have an obvious inclusion $\fhd\subset\fh_T^\star$. Hence we can define the element
$$
h_\lambda:=(\tau,\lambda)_T\in T
$$
for any $\lambda\in\fhd$. 

\subsection{Deformed category $\CO$}
Let $T$ be a local deformation algebra. 
For any $U\otimes_\DC T$--module $M$ and $\lambda\in\fhd$ set 
$$
M_\lambda:=\{m\in M\mid H.m=(\lambda+\tau)(H)m\quad\forall H\in\fh\},
$$
where we view $(\lambda+\tau)(H)$ as an element in $T$.  Let  $\CO_T$ be the category of $U\otimes_\DC T$--modules $M$ such that 
 $M=\bigoplus_{\lambda\in\fhd}M_\lambda$ and such that $B\otimes_\DC T.m$ is finitely generated over $T$ for any $m\in M$.
Then $\CO_T$ is an abelian category and in the special case $T=\DC$ with any $S$--structure we arrive at the well--studied BGG--category $\CO$. For any morphism of local deformation algebras $T\to T^\prime$ there is a base change functor
$$
\cdot\otimes_T T^\prime\colon\CO_{T}\to\CO_{T^\prime}.
$$
 For any $\lambda\in\fhd$ we define the Verma module
$$
M_T(\lambda):=U\otimes_B T_\lambda,
$$
where $T_\lambda$ denotes the $B$--structure on $T$ given by the composition $B\to S\stackrel{\lambda+\tau}\to T$ with  a left invers map $B\to S$  of the inclusion $S\to B$. The Verma modules are objects of $\CO_T$. 

For any $\nu\in\fhd$ consider the full subcategory   $\CO^{\varle\nu}_{T}$ of $\CO_T$ of modules $M$ such that $M_\mu=0$ if $\mu\not\le\nu$. In contrast to $\CO_T$, the  subcategories  $\CO^{\varle\nu}_{T}$ have enough projective objects. Moreover, the categories $\CO^{\varle\nu}_{T}$ provide a filtration of $\CO_T$ in the sense that every finitely generated object of $\CO_T$ lies in a finite direct sum of $\CO_T^{\varle\nu}$'s. 
Let $\fm\subset T$ be the maximal ideal of $T$ and $\DK=T/\fm$ the residue field. Consider the base change functor $\cdot\otimes_T\DK$.
\begin{theorem}[\cite{Fie}, Proposition 2.1 \& 2.6]
\indent
\begin{enumerate} 
\item The base change $\cdot\otimes_T \DK$ gives a bijection
$$
\left\{
\begin{matrix}
\text{simple isomorphism } \\
\text{classes of $\CO_T$} 
\end{matrix}
\right\}
\longleftrightarrow
\left\{
\begin{matrix}
\text{simple isomorphism } \\
\text{classes of $\CO_{\DK}$}
\end{matrix}
\right\}.
$$

\item For any $\nu\in\fhd$ the base change $\cdot\otimes_T \DK$ gives a bijection
$$
\left\{
\begin{matrix}
\text{projective isomorphism} \\
\text{ classes of $\CO_T^{\varle\nu}$} 
\end{matrix}
\right\}
\longleftrightarrow
\left\{
\begin{matrix}
\text{projective isomorphism} \\
\text{ classes of $\CO_{\DK}^{\varle\nu}$}
\end{matrix}
\right\}.
$$
\end{enumerate}
\end{theorem}
The category $\CO_\DK$ is a direct summand of the usual category $\CO$ over the Kac-Moody algebra $\fg\otimes_\DC \DK$ consisting of all objects whose weights lie
in the complex affine subspace  $\tau+\fhd=\tau+\Hom_\DC(\fh,\DC)\subset\fh_\DK^\star=\Hom_\DK(\fh\otimes_\DC \DK,\DK)$. Hence the simple isomorphism classes of $\CO_\DK$ (and hence of $\CO_T$) are parametrized by their highest weights, i.e.ø by elements of $\fhd$. Let $L_T(\lambda)$ be a simple object in $\CO_T$ corresponding to $\lambda\in\fhd$. It is a quotient of the Verma module $M_T(\lambda)$. 

We have the following structure theorem for projective objects.
\begin{theorem}[\cite{Fie}, Proposition 2.4 \& 2.7]\label{structure of projectives}
Let $T$ be a local deformation algebra and $\DK$ its residue field. Let $\nu\in\fhd$  and let $L_T(\lambda)$ be a simple object in $\CO_T^{\varle\nu}$.
\begin{enumerate}
\item
 There is a projective cover $P_T^{\varle\nu}(\lambda)$  of  $L_T(\lambda)$ in $\CO_T^{\varle\nu}$. Every projective object in $\CO_T^{\varle\nu}$ is isomorphic to a direct sum of projective covers. 
\item $P_T^{\varle\nu}(\lambda)$ has Verma flag, i.e.ø a filtration with subquotients isomorphic to Verma modules, and for the multiplicities the BGG--reciprocity formula  
$$
\left( P_T^{\varle\nu}(\lambda):{M_T}(\mu)\right) =\bigl[
{M_\DK}(\mu):L_\DK(\lambda) \bigr] 
$$
holds for  all Verma modules ${M_T}(\mu)$ in $\CO_T^{\varle\nu}$.
\item\label{structure of projectives - coherence}
Let $T\to T^\prime$ be a morphism of local deformation algebras and  $P$ projective  in $\CO_T^{\varle\nu}$. Then $P\otimes_T T^\prime$ is projective in $\CO_{T^\prime}^{\varle\nu}$. If $P$ is finitely generated,  then the natural transformation
$$
\Hom_{\CO_T}(P,\cdot)\otimes_T T^\prime \to \Hom_{\CO_{T^\prime}}(P\otimes_T T^\prime, \cdot\otimes_T T^\prime)
$$
is an isomorphism of functors from $\CO_T$ to $T^\prime\catmod$. 
\end{enumerate}
\end{theorem}
Part (\ref{structure of projectives - coherence}) of the theorem says that this deformation theory is a coherent deformation of categorical structures.
\subsection{Block decomposition}
Let $T$ be a local deformation algebra and $\DK$ its residue field.
Let $\sim_T$ be the equivalence relation on $\fhd$ generated by $\lambda\sim_T \mu$ if there exist $n\in \DN$ and $\beta\in\Delta_+$ such that 
$2(\lambda+\tau+\rho,\beta)_\DK=n(\beta,\beta)_\DK$ and $\lambda-\mu=n\beta$, where $\rho\in\fhd$ is a Weyl vector, i.e.ø $(\rho,\alpha)=1$ for any simple root $\alpha\in\Pi$. Then $\sim_T$ does not depend on the choice of $\rho$ and, by definition, $\sim_T=\sim_\DK$. For any union of equivalence classes $\Lambda\subset\fhd/_{\textstyle{\sim_T}}$ let $\CO_{T,\Lambda}$ be the full subcategory of $\CO_T$ consisting of all $M$ such that the highest weight of every simple subquotient of $M$ lies in $\Lambda$. If $\Lambda$ is a single equivalence class, then  $\CO_{T,\Lambda}$ is called a block of $\CO_T$.  

\begin{theorem}[\cite{Fie}, Proposition 2.8]
The functor $\{M_\Lambda\}  \mapsto  \bigoplus M_\Lambda $ is an equivalence of categories
$$
\prod_{\Lambda\in\fhd/_{\sim_T}}\CO_{T,\Lambda}\stackrel{\sim}\to  \CO_T.
$$
\end{theorem}

For any morphism of local deformation algebras $T\to T^\prime$ the equivalence relation $\sim_{T^\prime}$ is finer than $\sim_T$ and hence the base change $\cdot \otimes_T T^\prime$ respects the block decomposition, i.e.ø it induces a  base change functor $\CO_{T,\Lambda}\to\CO_{T^\prime,\Lambda}$.

We will consider only blocks outside the critical hyperplanes, i.e.ø blocks corresponding to equivalence classes which do not intersect the hyperplanes defined by $2(\lambda+\rho+\tau,\beta)_\DK=n(\beta,\beta)_\DK$ for any $n\in\DN$ and an {\em imaginary} root $\beta\in\Delta^{im}$. For any equivalence class
$\Lambda\in\fhd/_{\textstyle{\sim_T}}$ let
$$
\Delta_T(\Lambda):=\{\beta\in\Delta\mid 2(\lambda+\rho+\tau,\beta)_\DK\in\DZ(\beta,\beta)_\DK\quad\text{for some $\lambda\in\Lambda$}\}
$$
be the set of integral roots with respect to $\Lambda$. 
Hence $\Lambda$ lies outside the critical hyperplanes if and only if $\Delta_T(\Lambda)\subset \Delta^{re}$. In this case choose $\lambda\in\Lambda$. Let $\CW$ be the Weyl group of $\fg$. It naturally acts on $\fhd$ by linear transformations. Define the shifted action  by $w.\nu=w(\nu+\rho)-\rho$ for $w\in\CW$ and $\nu\in\fhd$. Again this action does not depend on the choice of $\rho$. 
 Then 
$$
\Lambda=\CW_T(\Lambda).\lambda,
$$
where $\CW_T(\Lambda)$ is the integral Weyl group with respect to  $\Lambda$, i.e.ø the group that is generated by the reflections $s_\alpha$ corresponding to  $\alpha\in\Delta_T(\Lambda)$. 

View $\DC$ as the residue field of $R$, hence it inherits the $S$--structure given by evaluation at $0\in\fhd$. Then $\CO_\DC$ is nothing else than the well--known  BGG--category $\CO$ over $\fg$ with respect to $\fh$ and $\fb$, $M_\DC(\lambda)$ and $L_\DC(\lambda)$ are  the Verma module and the simple module with highest weight $\lambda$ and the relations $\sim_\DC$ and $\sim_R$ coincide with the usual equivalence relation $\sim$ on $\fhd$, given by the submodule structure of Verma modules. Hence we will omit the subscripts and  write $\Delta(\Lambda)$ and $\CW(\Lambda)$ for the integral roots and the integral Weyl group of an equivalence class $\Lambda$ with respect to this relation. Define $\Delta_+(\Lambda)=\Delta(\Lambda)\cap\Delta_+$, $\Pi(\Lambda):=\{\alpha\in\Delta_+(\Lambda)\mid s_\alpha(\Delta_+(\Lambda)\setminus\{\alpha\})\subset\Delta_+(\Lambda)\}$ and $\CS(\Lambda)=\{s_\alpha\mid\alpha\in\Pi(\Lambda)\}$. Then $(\CW(\Lambda),\CS(\Lambda))$ is a Coxeter system. The main result of this article is that the structure of $\CO_\Lambda$ only depends on the isomorphism class of this Coxeter system and its action on $\Lambda$. 

\subsection{Structure of generic and subgeneric blocks}\label{subsection - structure of generic and subgeneric blocks}
Let $T$ be a local deformation algebra and $\Lambda\in\fhd/{\textstyle{\sim_T}}$ an equivalence class. We will call $\Lambda$ {\em generic} if it contains only one element, and {\em subgeneric} if it contains exactly two elements. The corresponding blocks will also be called generic and subgeneric, resp.   
In this section we will explicitely describe the structure of generic and subgeneric blocks of $\CO_T$. 

Suppose   $\Lambda\in\fhd/_{\textstyle{\sim_T}}$ contains only finitely many elements. Then $\CO_{T,\Lambda}$ is the same as  $\CO_{T,\Lambda}^{\varle\nu}$ for $\nu$ big enough, hence there exist projective covers $P_T(\lambda)$ of all simple objects $L_T(\lambda)$ in $\CO_{T,\Lambda}$ and their direct sum $P_{T,\Lambda}=\bigoplus_{\lambda\in\Lambda} P_T(\lambda)$  is a faithful, small projective in the sense of \cite{Mit}. The functor 
$$
\Hom(P_{T,\Lambda},\cdot)\colon \CO_{T,\Lambda}\to \catmodr\End(P_{T,\Lambda})
$$
is an equivalence of categories. We will describe a generalization of this equivalence for infinite $\Lambda$ in section \ref{section - combinatorics of O}. 

Suppose  $\Lambda\in\fhd/_{\textstyle{\sim_T}}$ is generic, i.e.ø $\Lambda=\{\lambda\}$. Then $P_{T,\Lambda}=P_T(\lambda)=M_T(\lambda)$. 
If $\Lambda$ is subgeneric, i.e.~ $\Lambda=\{\lambda,\mu\}$ and $\lambda>\mu$, then $P_T(\lambda)=M_T(\lambda)$ and there exists a short exact sequence 
$$
 0\to M_T(\lambda)\stackrel{i}\longrightarrow P_T(\mu)\stackrel{\pi}\longrightarrow M_T(\mu)\to 0.
$$
Standard arguments (involving the Jantzen filtration) show that there is an injection $M_\DK(\mu)\to M_\DK(\lambda)$ (note that there is, in general, no non--trivial map $M_T(\mu)\to M_T(\lambda)$). Composing with the surjection $P_\DK(\mu)\to M_\DK(\mu)$ gives a map $P_\DK(\mu)\to M_\DK(\lambda)$. By Theorem \ref{structure of projectives}, (\ref{structure of projectives - coherence}) we can lift it to a morphism $P_T(\mu)\to M_T(\lambda)$. Let $h\in\End(M_T(\lambda))=T$ be the composition with the inclusion of the short exact sequence above. Note that $h\in T$ is well--defined up to multiplication with an invertible element in $T$. 

\begin{proposition}[\cite{Fie}, Example 2.2 \& Proposition 3.4] \label{structure of blocks} 
\indent
\begin{enumerate}
\item{(generic case)} If $\Lambda=\{\lambda\}$ , then
$$
\CO_{T,\Lambda}\cong T\catmod.
$$
\item{(subgeneric case)} If $\Lambda=\{\lambda,\mu\}$  and $\lambda>\mu$, then 
$\End(P_{T,\Lambda})$ is isomorphic to the path algebra of the quiver 
\vskip .2cm
\hskip 2.5cm\hbox{
 \includegraphics{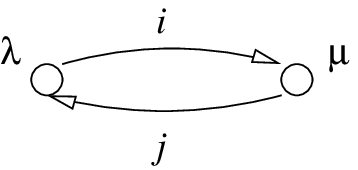} 
}

\noindent
over $T$ with relation  $j\circ i=h$, i.e.ø to the $T$--algebra generated by the paths of the quiver with relation $j\circ i=h e_\lambda$, where $e_\lambda$ is the trivial path at the vertex $\lambda$. 
\end{enumerate}
\end{proposition}
\begin{remark}\label{remarks on structure of subgeneric blocks} 
\begin{enumerate} 
\item In the case $T=R$ and $\lambda=s_\alpha.\mu$ for $\alpha\in\Delta^{re}$ it was shown in Proposition 3.4 in \cite{Fie}  that the homomorphisms can be chosen in such a way that $h=h_\alpha\in R$. 
\item  The quiver encodes the homomorphisms between $P_T(\mu)$ and $P_T(\lambda)$ as follows. The $\mu$--vertex (resp. the $\lambda$--vertex) in the quiver represents $P_T(\mu)$ (resp. $P_T(\lambda)$) and the space of morphisms is generated over $T$ by the paths between the resp. vertices modulo the relations. 
\end{enumerate} 
\end{remark}
\subsection{Translation functors}\label{subsection - translation functors}
Let $T$ be any local deformation algebra and $\Lambda, \Lambda^\prime\in\fhd/_{\textstyle{\sim_T}}$ equivalence classes outside the critical hyperplanes. Choose $\lambda\in\Lambda$ and $\lambda^\prime\in\Lambda^\prime$. Suppose that
\begin{enumerate}
\item $\lambda-\lambda^\prime$ is integral and either positive or negative. Hence $\Delta_T(\Lambda)=\Delta_T(\Lambda^\prime)$ and $\CW_T(\Lambda)=\CW_T(\Lambda^\prime)$.
\item Both $\lambda$ and $\lambda^\prime$ lie in the closure of the same Weyl chamber, i.e.ø $(\lambda+\rho+\tau,\beta)_\DK\in \DZ_{\ge 0}(\beta,\beta)$ if and only if $(\lambda^\prime+\rho+\tau,\beta)_\DK\in \DZ_{\ge 0}(\beta,\beta)$ for all $\beta\in\Delta_T(\Lambda)$.
\item $\Stab(\lambda)\subset \Stab(\lambda^\prime)$ with finite index ($\Stab$ denotes the stabilizer under the $\rho$--shifted action of $\CW$). Hence $\Lambda^\prime$ lies ``on the walls".
\end{enumerate} 
Let $\CM_T$ be  the full subcategory of $\CO_T$ consisting of modules which admit a finite Verma flag, i.e.~ a finite filtration whose subquotients are isomorphic to Verma modules. Let $\CM_{T,\Lambda}$ and $\CM_{T,\Lambda^\prime}$ be the corresponding blocks. In \cite{Fie} we defined the {\em translation functors}
$$
\vartheta_{out}=\vartheta_{T,out}\colon \CM_{T,\Lambda^\prime}Ê \to  \CM_{T,\Lambda} 
$$
and
$$
\vartheta_{on}=\vartheta_{T,on}\colon \CM_{T,\Lambda}Ê\to  \CM_{T,\Lambda^\prime}.
$$
In the case that $\Lambda$ is regular, i.e. $\Stab(\lambda)=\{e\}$ and $\Lambda^\prime$ is subregular, i.e.~ $\Stab(\lambda^\prime)=\{e,s\}$ for some reflection $s$, we call the composition $\Theta_s:=\vartheta_{out}\circ\vartheta_{on}$ the ``translation through the $s$--wall''.
We summarize the properties of translation functors in the following

\begin{theorem}[\cite{Fie}, Proposition 4.1 \& 5.2, Lemma 5.4,  Corollary 5.10 \& 5.11] \label{properties of translation functors}
\indent
\begin{enumerate} 
\item\label{properties of translation functors - exactness} $\vartheta_{out}$ and $\vartheta_{on}$ transform short exact sequences to short exact sequences.
\item \label{properties of translation functors - weights}
Let $w\in\CW_T(\Lambda)$. Then $\vartheta_{on} M_T(w.\lambda)$ is isomorphic to $M_T(w.\lambda^\prime)$ and $\vartheta_{out} M_T(w.\lambda^\prime)$ has a Verma flag with subquotients $M_T(w\bar x.\lambda)$ for $\bar x\in \Stab(\lambda^\prime)/\Stab(\lambda)$, each occurring once. 
\item If $\Stab(\lambda)=\Stab(\lambda^\prime)$, then $\vartheta_{on}$ and $\vartheta_{out}$ are mutually inverse equivalences of categories. 
\item Any base change $T\to T^\prime$ induces isomorphisms 
\begin{eqnarray*}
\vartheta_{T,on}(\cdot)\otimes_T T^\prime & \cong & \vartheta_{T^\prime,on}(\cdot\otimes_T T^\prime),
\\
\vartheta_{T,out}(\cdot)\otimes_T T^\prime & \cong & \vartheta_{T^\prime,out}(\cdot\otimes_T T^\prime)
\end{eqnarray*}
between the resp. translation functors (note that, in general, the functor $\vartheta_{T^\prime,\ast}(\cdot\otimes_T T^\prime)$ splits into a direct product of translation functors, cp. \cite{Fie}, Lemma 5.4).
\item $\vartheta_{on}$ is left adjoint to $\vartheta_{out}$.  If $P\in\CM_{T,\Lambda}$ is projective in $\CO_{T,\Lambda}$, then $\vartheta_{on} P$ is projective in $\CO_{T,\Lambda^\prime}$.
\item\label{properties of translation functors - biadjointness} Assume that $T$ is $R$, $R_\fp$ or one of their residue fields. Then $\vartheta_{out}$ and $\vartheta_{on}$ are biadjoint. If $P^\prime\in\CM_{T,\Lambda^\prime}$ is projective in $\CO_{T,\Lambda^\prime}$, then $\vartheta_{out} P^\prime$ is projective in $\CO_{T,\Lambda}$. 
\end{enumerate}
\end{theorem}

\subsection{The tilting equivalence}\label{subsection - tilting equivalence}
Let $T$ be any local deformation algebra with $S$--structure $\tau\colon S\to T$. Let $\gamma\colon S\to S$ be the automorphism which is given by $\gamma(h)=-h$ for all $h\in\fh$ and let $\bar T$ be the $S$--algebra which is isomorphic to $T$ as an algebra, but where the $S$--structure is replaced by $\bar\tau=\tau\circ\gamma$. In particular, $\gamma$ can be thought of as an isomorphism $\gamma\colon S\to\bar S$ of $S$--algebras. Since $\fh\subset S$ is $\gamma$--stable, it induces an isomorphism of $S$--algebras $R\to\bar R$ and for any $\gamma$--stable prime ideal $\fp\subset R$ an isomorphism $R_\fp\to\bar R_\fp$. These are again denoted by $\gamma$.

Following \cite{Soe2} we will construct the tilting functor  $t^\prime\colon\CM_T\to\CM_{\bar T}^{opp}$ corresponding to the  semi--infinite character $2\rho$. Let $S_{2\rho}$ be the semi--regular $U$--bimodule of \cite{Soe2}. For any $\fg\otimes T$--module $N$ we define $N^\circledast:=\bigoplus_{\lambda\in\fhd}\Hom_T(N_\lambda,T)$ and we let $\fg$ act on $N^\circledast$ by defining $(Xf)(n)=-f(Xn)$ for $n\in N$, $f\in N^\circledast$ and $X\in\fg$.
For any $M\in\CM_T$ set 
$$
t^\prime(M):=(S_{2\rho}\otimes_U M)^\circledast.
$$ 
As in \cite{Soe2} one shows that $t^\prime(M_T(\lambda))$ is isomorphic to $M_{\bar T}(-2\rho-\lambda)$, that $t^\prime$ transforms short exact sequences of modules with Verma flag to short exact sequences and that it defines an equivalence  $t^\prime\colon\CM_T\to\CM_{\bar T}^{opp}$ of categories.

Let $T$ be $R$, $R_\fp$, where $\fp$ is $\gamma$--stable, or one of their residue fields. 
The isomorphism $\gamma\colon T\to\bar T$  induces an equivalence of categories $\CM_T\cong \CM_{\bar T}$ and composition with $t^\prime$ gives an equivalence $t\colon \CM_T\to\CM_T^{opp}$ with $t(M_T(\lambda))\cong M_T(-2\rho-\lambda)$. Moreover, $t$ respects the block structure, i.e.~ for any $\Lambda\in\fhd/_{\textstyle{\sim_T}}$ it induces an equivalence 
$$
t=t_\Lambda\colon\CM_{T,\Lambda}\stackrel{\sim}\to\CM_{T,t(\Lambda)},
$$
where $t(\Lambda):=\left\{-2\rho-\lambda\mid \lambda\in\Lambda\right\}$. An easy calculation shows that this is again an equivalence class (note that $-w.\lambda-2\rho=w.(-\lambda-2\rho)$). This is most useful for us, since $t(\Lambda)$ contains a highest weight if and only if $\Lambda$ contains a lowest weight. This allows us to choose for any of the following constructions the most convenient case. We will, for example, first define the structure functor for equivalence classes with lowest weight, but we will then use it mainly in the case with highest weight, because this is the case where projective objects exist.

Let $Z_{T,\Lambda}$ be the center of $\CO_{T,\Lambda}$, i.e.~ the ring of endotransformations of the identity functor. In Section \ref{subsection - generic and subgeneric cases} we will give an explicit description of this center for generic and subgeneric blocks and in Theorem \ref{calculation of the center} for any block outside the critical hyperplanes provided that $T=R$.
\begin{lemma}\label{tilted center}
There is a natural identification
$$
Z_{T,\Lambda}\stackrel{\sim}\to Z_{T,t(\Lambda)}
$$
induced by $t$.
\end{lemma}
\begin{proof}
Restriction gives an isomorphism between the center of $\CO_{T}$ and the center of $\CM_T$, since every finitely generated object of  $\CO_{T}$ admits a resolution by objects with finite Verma flag. Moreover, the center of any category is the center of its opposed category.
\end{proof}
\begin{remark} In the case $T=R$ 
Theorem \ref{calculation of the center} gives a simple and  explicit construction of this identification. Let $\lambda\in\Lambda$ and $z=\{z_w\}\in Z_{R,\Lambda}$, so $z$ acts on $M_R(w.\lambda)$ as multiplication with $z_w$. Then the image of $z$ under the identification in the lemma will act as multiplication with $z_w$ on $M_R(-w.\lambda-2\rho)$.
\end{remark}

\section{The Struktursatz and the combinatorics of translation functors}

\subsection{The generic and subgeneric cases}\label{subsection - generic and subgeneric cases}
Let $T$ be a local deformation algebra.  If  $\Lambda=\{\lambda\}$ is generic, then $P_T(\lambda)=M_T(\lambda)$ and the evaluation $Z_{T,\Lambda}\to\End(P_T(\lambda))=T$ is an isomorphism. The functor 
$$
\DV=\DV_{T,\Lambda}:=\Hom(P_T(\lambda),\cdot)\colon\CO_{T,\Lambda}\to T\catmod
$$
is an equivalence of categories. 

Suppose $\Lambda=\{\lambda,\mu\}$ is subgeneric and $\lambda>\mu$. In Section \ref{subsection - structure of generic and subgeneric blocks} we defined an element $h\in T$.
\begin{lemma} \label{calculation of the center - subgeneric case}
\begin{enumerate}
\item The evaluation
 $Z_{T,\Lambda}\stackrel{\sim}\to\End(P_T(\mu))$  at the antidominant projective object is an isomorphism. 
\item Suppose $h\ne 0$. Then the evaluation
$$
Z_{T,\Lambda}\to\End(M_T(\lambda))\oplus\End(M_T(\mu))=T\oplus T
$$
is injective and its image is $\left\{(t_\lambda,t_\mu)\mid t_\lambda\equiv t_\mu\mod h\right\}$.
\end{enumerate}
\end{lemma}
\begin{proof} The explicit description of $\End(P_T(\mu))$ in Proposition \ref{structure of blocks} shows that it is generated over $T$ by the identity and the composition $i\circ j$ with relation $(i\circ j)^2=h(i\circ j)$. In particular it is commutative.
 
Let $f\in\End(P_T(\mu))$. Consider the short exact sequence 
$$
 0\to M_T(\lambda)\to P_T(\mu)\to M_T(\mu)\to 0.
$$
Then $f(M_T(\lambda))\subset M_T(\lambda)$ 
since $\lambda>\mu$, hence $f$ induces a map $f_\lambda\in\End(M_T(\lambda))$ and a map $f_\mu\in\End(M_T(\mu))$. Since any element of the center is determined by its actions on the indecomposable projective objects $P_T(\mu)$ and $M_T(\lambda)$, we conclude that the map $Z_{T,\Lambda}\to\End(P_T(\mu))$ is injective. Conversely, a pair $(f,g)\in\End(P_T(\mu))\oplus \End(M_T(\lambda))$ which commutes with homomorphisms defines an element of the center. For any $f\in\End(P_T(\mu))$ the pair $(f,f_\lambda)$ is of such type, hence the evaluation map is surjective and we proved the first statement. 

The evaluation in the second statement clearly factors over the map 
\begin{eqnarray*}
\End(P_T(\mu)) &  \to & \End(M_T(\lambda))\oplus\End(M_T(\mu))=T\oplus T\\
f & \mapsto & (f_\lambda,f_\mu)
\end{eqnarray*}
that we just defined. 
It sends the identity to $(1,1)$ and $i\circ j$ to $(h,0)$. Hence, 
if $h\ne 0$, this map is injective and its image is the set  of pairs $(t_\lambda,t_\mu)$ such that $t_\lambda\equiv t_\mu$ modulo $h$.
\end{proof}
Define the {\em structure functor}
$$
\DV=\DV_{T,\Lambda}:=\Hom(P_T(\mu),\cdot)\colon \CO_{T,\Lambda}\to Z_{T,\Lambda}\catmod.
$$
The next theorem is the subgeneric case of Soergel's Struktursatz. Note, however, that we do not assume that $M$ or $M^\prime$ is projective.
\begin{theorem}\label{struktursatz - subgeneric case} Suppose that $\Lambda$ is subgeneric and that $h\ne0$. Then $\DV$ is fully faithful on $\CM_{T,\Lambda}$, i.e.~ for any $M,M^\prime\in\CM_{T,\Lambda}$ the map
$$
\Hom_{\CM}(M,M^\prime)\stackrel{\sim}\to\Hom_{Z}(\DV M,\DV M^\prime)
$$
is an isomorphism. 
\end{theorem}
\begin{proof} From the explicit description of the structure of $\CO_{T,\Lambda}$ in Proposition \ref{structure of blocks} follows that $M_T(\lambda)$, $M_T(\mu)$ and $P_T(\mu)$ are the only indecomposable objects in $\CM_{T,\Lambda}$. 
We will now define their counterparts in $Z_{T,\Lambda}\catmod$.
Let $\FM(\lambda)=T$ and $\FM(\mu)=T$ with the $Z_{T,\Lambda}$--action given by letting $(z_\lambda,z_\mu)\in Z_{T,\Lambda}$ act by multiplication with $z_\lambda$ and $z_\mu$, resp. Let $\FP$ be $Z_{T,\Lambda}$, considered as a $Z_{T,\Lambda}$--module. Then $\DV(M_T(\lambda))\cong \FM(\lambda)$, $\DV(M_T(\mu))\cong \FM(\mu)$ and $\DV(P_T(\mu))\cong \FP$. Define maps
$$
\begin{matrix}
a\colon \FM(\lambda) \to\FP, &  1 \mapsto (h,0) \\
b\colon \FP\to\FM(\lambda), & (z_\lambda, z_\mu)  \mapsto z_\lambda \\
c\colon \FM(\mu) \to\FP, &  1 \mapsto (0,h) \\
d\colon \FP\to\FM(\mu), & (z_\lambda, z_\mu)  \mapsto z_\mu. \\
\end{matrix}
$$
It is easy to check that these maps generate all homomorphisms and that a complete set of relations  is given by the following: 
$d\circ a=0$, $b\circ c=0$ and $b\circ a$, $d\circ c$ and $a\circ b+c\circ d$ are multiplication with $h$
 (note, in particular, that there are no non--trivial homomorphism between $\FM(\lambda)$ and $\FM(\mu)$).

The elements $(h,0)$ and $(0,h)$ in $Z_{T,\Lambda}$ acts on $P_T(\mu)$ and provide maps $b^\prime\colon P_T(\mu)\to M_T(\lambda)$ and $c^\prime\colon M_T(\mu)\to P_T(\mu)$.  Let $a^\prime=i\colon M_T(\lambda)\to P_T(\mu)$ and $d^\prime=\pi\colon P_T(\mu)\to M_T(\mu)$ be the injection and surjection of the short exact sequence in Section \ref{subsection - structure of generic and subgeneric blocks}. Using Proposition \ref{structure of blocks} again one checks that these maps generate all homomorphisms and fulfill the same relations as before.
\end{proof}

\begin{remark} In the non--deformed situation, i.e.~ in the case $T=\DC$, we unfortunately have $h=0$. For finite dimensional $\fg$ it was proven in \cite{Soe1} that the map above is still a bijection provided that $M^\prime$ is projective.
\end{remark}

Let $T$ be a local deformation algebra. Choose $\Lambda$ and $\Lambda^\prime$ as in Section \ref{subsection - translation functors}. If $\Lambda$ and $\Lambda^\prime$ are both generic or both subgeneric, then the translation functors $\vartheta_{on}$ and $\vartheta_{out}$ are mutually inverse equivalences of categories. Suppose now that $\Lambda=\{\lambda,\mu\}$ is subgeneric with $\lambda<\mu$ and that $\Lambda^\prime=\{\lambda^\prime\}$ is generic. The arguments applied to the case $T=R_\alpha$ in \cite{Fie} apply in fact to any subgeneric situation, i.e. all the statements in Theorem \ref{properties of translation functors} hold for arbitrary $T$ provided that $\Lambda$ is subgeneric. In particular we deduce that $\vartheta_{out}(P_T(\lambda^\prime))\cong P_T(\lambda)$.  The corresponding map $\End(P_T(\lambda^\prime))\to\End(P_T(\lambda))$ can be interpreted, using Lemma \ref{calculation of the center - subgeneric case}, as a map $Z_{T,\Lambda^\prime}\to Z_{T,\Lambda}$. This map is canonical, even though the identification $\vartheta_{out}(P_T(\lambda^\prime))\cong P_T(\lambda)$ is not. It is nothing else than the canonical inclusion $T\cdot 1\hookrightarrow Z_{T,\Lambda}$. Let $\Ind= \Ind_{Z_{T,\Lambda^\prime}}^{Z_{T,\Lambda}}=Z_{T,\Lambda}\otimes_{Z_{T,\Lambda^\prime}}\cdot$ be the induction functor and $\Res=\Res^{Z_{T,\Lambda^\prime}}_{Z_{T,\Lambda}}$ its right adjoint.

\begin{theorem}\label{comb of translation - subgeneric} Suppose $\Lambda$ is subgeneric and $\Lambda^\prime$ is generic. Then there are isomorphisms of functors
$$
\DV\circ\vartheta_{on}\cong\Res\circ\DV\colon \CM_{T,\Lambda}\to Z_{T,\Lambda^\prime}\catmod
$$
and
$$
\DV\circ\vartheta_{out}\cong\Ind\circ\DV\colon \CM_{T,\Lambda^\prime}\to Z_{T,\Lambda}\catmod.
$$
\end{theorem}
\begin{proof} We have
\begin{eqnarray*}
\DV\circ\vartheta_{on}M & = & \Hom(P_T(\lambda^\prime), \vartheta_{on}M) \\
& \cong & \Hom(\vartheta_{out} P_T(\lambda^\prime),  M) \\
& \cong & \Hom(P_T(\lambda),  M) \\
& = & \Res\circ\DV M
\end{eqnarray*}
and hence proved the first statement. It provides, together with the adjunctions $\id\to\vartheta_{on}\circ\vartheta_{out}$ and $\Ind\circ\Res\to\id$, a morphism
$$
\Ind\circ\DV\to\Ind\circ\DV\circ\vartheta_{on}\circ\vartheta_{out}\to\Ind\circ\Res\circ\DV\circ\vartheta_{out}\to\DV\circ\vartheta_{out}.
$$
We want to prove that this is an isomorphism.  It is enough to check that it is an isomorphism when evaluated at $M=M_{T}(\lambda^\prime)=P_{T}(\lambda^\prime)$. In this case the map $\DV P_T(\lambda^\prime)\to\Res\,\DV\vartheta_{out} P_T(\lambda^\prime)$ is, up to an invertible scalar, just the previously defined map $Z_{T,\Lambda^\prime}\to Z_{T,\Lambda}$ (of $Z_{T,\Lambda^\prime}$--modules).
After applying the induction $Z_{T,\Lambda}\otimes_{Z_{T,\Lambda^\prime}}\cdot$  and composing with the multiplication $Z_{T,\Lambda}\otimes_{Z_{T,\Lambda^\prime}} Z_{T,\Lambda}\to Z_{T,\Lambda}$ we get an isomorphism. Hence $\Ind\circ\DV\cong\DV\circ\vartheta_{out}$.
\end{proof}

\subsection{Splitting of equivalence classes under localizations}
Now we assume that $T=R$, i.e.~ the localization of $S=S(\fh)$ at the maximal ideal generated by $\fh\subset S$. Note that $\sim_R=\sim_\DC=\sim$ is the usual Kac-Kazhdan relation, so in the following we suppress the index $R$. We want to prove the statements in Theorems \ref{struktursatz - subgeneric case} and \ref{comb of translation - subgeneric} for arbitrary $\Lambda\in\fhd/_{\textstyle{\sim}}$ outside the critical hyperplanes and $T=R$. We will employ the base change functors $\cdot\otimes_R R_\fp$ for prime ideals $\fp\subset R$ of height one and use the coherence in Theorem \ref{structure of projectives} together with the facts that all objects we will work with are free over $R$ and that all constructions commute with base change. Most importantly, the following lemma shows that the localization $\cdot\otimes_R R_\fp$ reduces our situation to generic and subgeneric situations. Its proof can be found in \cite{Fie}, section 3.3.

\begin{lemma} \label{splitting of equivalence classes} Let $\Lambda\in\fhd/_{\textstyle{\sim}}$ (not necessarily outside the critical hyperplanes) and let $\fp\subset R$  be a prime ideal.
\begin{enumerate}
\item\label{equivalence classes 1}  If $h_\alpha\not\in\fp$ for all roots $\alpha\in\Delta(\Lambda)$, then $\Lambda$ splits under $\sim_{R_\fp}$ into trivial equivalence classes, i.e.ø into equivalence classes with only one element. 
\item \label{equivalence classes 2} If $\fp=R h_\alpha$ for a real root $\alpha\in\Delta(\Lambda)$, then $\Lambda$ splits under $\sim_{R_\fp}$ into equivalence classes of the form $\{\lambda, s_\alpha.\lambda\}$.
\end{enumerate}
\end{lemma}
If $\Lambda$ lies outside the critical hyperplanes, i.e.~ $\Delta(\Lambda)\subset\Delta^{re}$, and $\fp\subset R$ is a prime ideal of height one, either (1) or (2) applies. 

\subsection{The center of $\CO_{R,\Lambda}$}\label{subsection - the center} Let $T=R$ and 
$\Lambda\in\fhd/_{\textstyle{\sim}}$ be an equivalence class outside the critical hyperplanes and choose $\lambda\in\Lambda$, hence $\Lambda=\CW(\Lambda).\lambda$. 
In Section \ref{subsection - tilting equivalence} we defined the notion of the center $Z_{R,\Lambda}$ of $\CO_{R,\Lambda}$. Let 
$$
Z_{R,\Lambda}\to\prod_{w\in\CW(\Lambda)/\Stab(\lambda)}\End(M_R(w.\lambda))=\prod_{w\in\CW(\Lambda)/\Stab(\lambda)}R
$$
be the evaluation at the Verma modules. Remember the element $h_\alpha=(\tau,\alpha)_R\in R$.
\begin{theorem}[\cite{Fie}, Theorem 3.6] \label{calculation of the center} 
 The evaluation map is injective and induces an isomorphism
$$
Z_{R,\Lambda}\cong\left\{
{
\left.
(z_w)\in\prod_{w\in\CW(\Lambda)/\Stab(\lambda)}R \;\right\lvert 
\begin{matrix}
z_w\equiv z_{s_\alpha w}\mod h_\alpha \\
 \;\forall \alpha\in\Delta(\Lambda), w\in\CW(\Lambda)/\Stab(\lambda)
\end{matrix}
}
\right\}. 
$$
\end{theorem}
\begin{remark} In \cite{Fie}  it was shown that, after identifying the centers with their images under the evaluation map, we have $Z_{R,\Lambda}=\bigcap_\fp Z_{R_\fp,\Lambda}$, where the intersection is taken over all prime ideals of height one. Hence the theorem above follows from Lemma \ref{calculation of the center - subgeneric case} and Lemma \ref{splitting of equivalence classes} (cp. also Remark \ref{remarks on structure of subgeneric blocks}).
\end{remark}
Now suppose that $\Lambda\in\fhd/_{\textstyle{\sim}}$ lies outside the critical hyperplanes and, in addition, is regular, i.e.~ $\Stab(\lambda)=\{e\}$ for some (hence any) $\lambda\in\Lambda$. 
We then define a left action of $\CW(\Lambda)$ on $Z_{R,\Lambda}$. For
 $x\in\CW(\Lambda)$ and $z=(z_w)\in Z_{R,\Lambda}\subset\prod R$ define $x.z=(z^\prime_w)$ with $z^\prime_w=z_{wx}$.  For $\alpha\in\Delta(\Lambda)$ let $Z^{s_\alpha}_{R,\Lambda}\subset Z_{R,\Lambda}$ be the $s_\alpha$--invariants. 

\subsection{The structure functor}
In this subsection $T$ can be an arbitrary local deformation algebra.
 We will define a functor $\DV_{T,\Lambda}\colon\CM_{T,\Lambda}\to Z_{T,\Lambda}\catmod$. 

Let $\Lambda\in\fhd/_{\textstyle{\sim_T}}$ be an equivalence class outside the critical hyperplanes. We will call $\Lambda$ of  {\em negative level} if it contains an antidominant element, i.e.~ an element $\lambda$ with the property $\lambda\le\lambda^\prime$ for any $\lambda^\prime\in\Lambda$. We call $\Lambda$ of {\em positive level} if it contains a dominant element, i.e.~ an element $\lambda$ with $\lambda^\prime\le\lambda$ for any $\lambda^\prime\in\Lambda$. This is consistent with the notions of positive and negative level in the case of affine Kac-Moody algebras, i.e.~ $\Lambda$ is of positive (resp. negative) level if all its elements are of positive (resp. negative) level.

Suppose now that $\Lambda$ is of negative level and $\lambda\in\Lambda$ is its smallest element. We will define a projective limit of the modules $P_T^{\varle\nu}(\lambda)$. Let $\chi=\sum \alpha$ be the sum of all simple roots of $\fg$. Hence, for all $\lambda^{\prime}\in\Lambda$, $\lambda^\prime\le \lambda+n\chi$ if $n$ is big enough. We can choose
for any $n\in \DN$ a surjective morphism $P_T^{\varle\lambda+(n+1)\chi}(\lambda)\sur P_T^{\varle\lambda+n\chi}(\lambda)$ and define the $\fg\otimes_\DC T$--module
$$
P^\infty_T(\lambda):=\varprojlim P_T^{\varle\lambda+n\chi}(\lambda).
$$
Then $P^\infty_T(\lambda)$ is an object of $\CO_{T,\Lambda}$ if and only if $\Lambda$ is finite. This is the case if and only if $\varprojlim P_T^{\varle\lambda+n\chi}(\lambda)$ stabilizes.

\begin{lemma}\label{antidominant multiplicities} Suppose $\Stab(\lambda)$ is finite. Then $P_T^\infty(\lambda)$ has a reversed Verma flag, i.e.ø a descending filtration whose subquotients are isomorphic to Verma modules, and for the multiplicities holds
$$
\bigl(P_T^\infty(\lambda):M_T(w.\lambda)\bigr)=1
$$
for all $w\in\CW_T(\Lambda)$. 
\end{lemma}

\begin{proof} 
First suppose that $\Lambda$ is regular, i.e. $\Stab(\lambda)=\{e\}$. Each $P^{\varle\lambda+n\chi}_T(\lambda)$ has a Verma flag by Theorem \ref{structure of projectives}.  
It suffices to show that $\left(P_T^{\varle \lambda+n\chi}(\lambda), M_T(w.\lambda)\right)=1$ if $n$ is such that $\lambda+n\chi\ge w.\lambda$. By BGG--reciprocity we have $\left(P_T^{\varle \lambda+n\chi}(\lambda), M_T(w.\lambda)\right)=\bigl [M_\DK(w.\lambda): L_\DK(\lambda)\bigr]$, where $\DK$ is the residue field of $T$. In particular the left hand side is independent of $n$ for $n\gg0$. Denote it by $m_w$. Standard arguments (involving the Jantzen sum formula) show that there is, up to scalars, a unique injection $M_\DK(sw.\lambda)\hookrightarrow M_\DK(w.\lambda)$ if $s\in\CS(\Lambda)$ is simple and $sw.\lambda\le w.\lambda$. By induction we get inclusions $M_\DK(\lambda)\hookrightarrow M_\DK(w.\lambda)$ for all $w\in\CW(\Lambda)$,  hence $m_w\ge 1$ for all $w\in\CW(\Lambda)$. From the Jantzen sum formula we also deduce $m_s=1$ for $s\in\CS(\Lambda)$. 

Choose $s\in\CS(\Lambda)$ and let $\Theta_s=\vartheta_{out}\vartheta_{on}$ be the translation through the $s$--wall. Then $\Theta_s P_T^{\varle \lambda+n\chi}(\lambda)$ is projective in a suitably truncated category and, by Theorem \ref{properties of translation functors}, the subquotient $M_T(\lambda)$ occurs with multiplicity $m_s+m_e=2$, hence it must contain two copies of $P_T^{\varle \lambda+n\chi}(\lambda)$ as direct summands. More generally, the subquotient $M_T(w.\lambda)$ occurs with multiplicity $m_w+m_{ws}$ and we deduce $m_w+m_{ws}\ge 2m_w$, hence $m_w\le m_{ws}$ for all $w\in\CW(\Lambda)$ and all simple reflections $s\in \CS(\Lambda)$. By induction on the length of $w$ we get $m_w\le 1$ for all $w\in\CW(\Lambda)$, hence $m_w=1$ and we proved the lemma for regular $\Lambda$.

If $\Lambda$ is not regular we can choose a regular $\Lambda^\prime$ which contains an antidominant $\lambda^\prime$ with the properties needed for the definition of translation functors. Then, by what we proved before and Theorem \ref{properties of translation functors},  $\left(\vartheta_{on} P_T^{\varle\lambda^\prime+n\chi}(\lambda^\prime):M_T(w.\lambda)\right)=\#\Stab(\lambda)$ for $n$ big enough. Moreover, $\vartheta_{on} P_T^{\varle\lambda^\prime+n\chi}(\lambda^\prime)$ is projective in a suitably truncated category, hence $\varprojlim \vartheta_{on} P_T^{\varle\lambda^\prime+n\chi}(\lambda^\prime)$ must contain $\#\Stab(\lambda)$ copies of $P^\infty_T(\lambda)$. We deduce $\left(P_T^\infty(\lambda): M_T(w.\lambda)\right)\le 1$.  The same arguments as in the regular case show $\left(P_T^\infty(\lambda): M_T(w.\lambda)\right)\ge 1$ .

\end{proof}

\begin{remarks}\label{remarks on multiplicity result}
\begin{enumerate}

\item \label{remarks on multiplicity result - splitting} It follows that for any base change $T\to T^\prime$, $P_T^\infty(\lambda)\otimes_T T^\prime$ splits into the direct sum of antidominant projective covers of $\CO_{T^\prime,\Lambda}$ and each occurs once, i.e.~ we have an isomorphism
$$
P_T^\infty(\lambda)\otimes_T T^\prime\cong\bigoplus_i P^\infty_{T^\prime}(\lambda_i),
$$
where $\Lambda=\dot\bigcup_i\Lambda_i$ is the splitting of $\Lambda$ with respect to $\sim_{T^\prime}$ and $\lambda_i\in\Lambda_i$ is the antidominant element. 
\item\label{remarks on multiplicity result - transferred antidominant projective} Suppose $T=R, R_\fp$ or one of its residue fields.  Let $\Lambda$ and $\Lambda^\prime$  be two equivalence classes as in Section \ref{subsection - translation functors} and suppose both contain antidominant elements $\lambda\in\Lambda$ and $\lambda^\prime\in\Lambda^\prime$. Then $\varprojlim\vartheta_{out}P_T^{\varle\lambda^\prime+n\chi}(\lambda^\prime)$ has a Verma flag and every Verma module in $\CM_{T,\Lambda}$ occurs once. Since its truncations are projective, we conclude that  $P_T^\infty(\lambda)\cong \varprojlim\vartheta_{out}P_T^{\varle\lambda^\prime+n\chi}(\lambda^\prime)$.
\item\label{remarks on multiplicity result - kernel of truncation} We can also deduce that the kernel of the chosen projection $P_T^{\varle\lambda+(n+1)\chi}(\lambda)\to P_T^{\varle\lambda+n\chi}(\lambda)$ is generated by all Verma subquotients with highest weight $\nu$ such that $\nu\not\le\lambda+n\chi$. 
\end{enumerate}
\end{remarks}

The following theorem is a consequence of the lemma and was proved for $\Lambda$ in positive level by Kashiwara-Tanisaki in \cite{KT4}. It is not used in the sequel. Let $\Lambda\in\fhd/_{\textstyle{\sim}}$ be outside the critical hyperplanes. Note that the theorem refers to the non--deformed situation, i.e.~ to the case $T=\DC$. There are no nontrivial homomorphisms $M_R(\lambda)\to M_R(\mu)$ if $\lambda\ne\mu$. 

\begin{theorem} \label{Verma embeddings} Let $\Lambda\in\fhd/_{\textstyle{\sim}}$ be an equivalence class outside the critical hyperplanes and $\lambda\in\Lambda$ dominant or antidominant. Let $w,w^\prime\in\CW(\Lambda)$ with $w.\lambda\le w^\prime.\lambda$. Then 
$$
\dim_\DC\Hom(M(w.\lambda),M(w^\prime.\lambda))=1.
$$
\end{theorem}
\begin{proof} Since any Verma module is free over $U(\fn^-)$, where $\fn^-$ is the subalgebra of $\fg$ generated by all negative weight spaces, and since $U(\fn^-)$ is torsion free,  every non--trivial homomorphism between Verma modules is injective. Moreover, as shown in the proof of Lemma \ref{antidominant multiplicities}, $\Hom(M(w.\lambda),M(w^\prime.\lambda))$ is at least one--dimensional for $w.\lambda\le w^\prime.\lambda$.

For antidominant $\lambda$ we deduce $[M(w.\lambda):L(\lambda)]=1$ from Lemma \ref{antidominant multiplicities} and BGG--re\-ci\-pro\-ci\-ty, hence $\dim_\DC\Hom(M(\lambda),M(w.\lambda))=1$. We get $\dim_\DC\Hom(M(w.\lambda),M(\lambda))=1$ for all dominant $\lambda$ using the tilting equivalence. Since all non--trivial morphisms between Verma modules are injective, we deduce $\dim_\DC\Hom(M(w.\lambda),M(w^\prime.\lambda))=1$ for $w,w^\prime\in\CW(\Lambda)$ such that $w.\lambda \le w^\prime.\lambda$.  Again using the tilting equivalence gives the statement for antidominant $\lambda$.
\end{proof}

We return to the assumption  that $\Lambda$ is of negative level. Let $\lambda\in\Lambda$ be the smallest element. The action of $Z_{T,\Lambda}$ on $P_T^{\varle\lambda+n\chi}(\lambda)$ induces a map $Z_{T,\Lambda}\to\End(P^\infty_T(\lambda))$. Hence we can define
\begin{eqnarray*}
 \DV=\DV_{T,\Lambda}\colon \CO_{T,\Lambda}  & \to &   Z_{T,\Lambda}\catmod \\
M & \mapsto & \Hom(P^\infty_T(\lambda), M).
\end{eqnarray*} 
For convenience we will always restrict $\DV$ to the subcategory $\CM_{T,\Lambda}$ of modules which admit a Verma flag and denote this restriction by the same symbol.

\begin{proposition}\label{properties of V}
\begin{enumerate}
\item There is an  isomorphism of functors
$$
\DV|_{\CM_{T,\Lambda}^{\varle\lambda+n\chi}}\cong\Hom(P_T^{\varle\lambda+n\chi}(\lambda), \cdot)\colon \CM_{T,\Lambda}^{\varle\lambda+n\chi}\to Z_{T,\Lambda}\catmod.
$$
In particular, $\DV$ transforms short exact sequences to short exact sequences.
\item  \label{properties of V - Vermas} For any $w\in\CW_T(\Lambda)$ there is an isomorphism $\DV M_T(w.\lambda)\cong Z_{T,\Lambda}/\fm_w$, where $\fm_w\subset Z_{T,\Lambda}$ is the ideal generated by all elements acting trivially on $M_T(w.\lambda)$. In particular, $\DV M_T(\lambda)$ is free of rank one  over $T$.
\item \label{properties of V - freeness over T} If $M\in\CM_{T,\Lambda}$, then $\DV M$ is free of finite rank over $T$.
\item \label{properties of V - base change} $\DV$ commutes with base changes $T\to T^\prime$, i.e.ø there are natural isomorphisms of functors 
$$
\DV_{T,\Lambda}(\cdot)\otimes_T T^\prime\cong\prod_{i}\DV_{T^\prime,\Lambda_i}(\cdot\otimes_T T^\prime),
$$
where $\Lambda=\dot\bigcup_i\Lambda_i$ is the splitting of $\Lambda$ under $\sim_{T^\prime}$.
\end{enumerate}
\end{proposition}
\begin{proof} 

The first statement directly follows Remark \ref{remarks on multiplicity result}, (\ref{remarks on multiplicity result - kernel of truncation}).

Let $\DK$ be the residue field of $T$. By  Theorem \ref{structure of projectives}, (\ref{structure of projectives - coherence})  and Lemma \ref{antidominant multiplicities} we have $\dim_\DK (\DV M_T(w.\lambda))\otimes_T\DK=\dim_\DK \Hom(P^{\infty}_\DK(\lambda), M_\DK(w.\lambda))=1$. Hence, by Nakayama's lemma, there is a generator $f$ of $\DV M_T(w.\lambda)$ as a $T$--module. Then $f$ is killed by the action of $\fm_w\subset Z_{T,\Lambda}$, hence the action of $Z_{T,\Lambda}$ on $f$ provides an isomorphism   $\DV M_T(w.\lambda)\cong Z_{T,\Lambda}/\fm_w=T$, hence we proved the second statement. 

The third statement follows immediately from the fact that $\DV$ transforms short exact sequences to short exact sequences and  (\ref{properties of V - Vermas}).

By Remark \ref{remarks on multiplicity result}, (\ref{remarks on multiplicity result - splitting}), the localization  $P_T^{\varle\lambda+n\chi}(\lambda)\otimes_T T^\prime$ splits into the direct sum of all antidominant projective covers in $\CO^{\varle\lambda+n\chi}_{T^\prime,\Lambda}$. By Theorem \ref{structure of projectives}, (\ref{structure of projectives - coherence}), the natural map defines an isomorphism of functors 
$$
\Hom(P_T^{\varle\lambda+n\chi}(\lambda),\cdot)\otimes_T T^\prime= \Hom(P_T^{\varle\lambda+n\chi}(\lambda)\otimes_T T^\prime,\cdot \otimes_T T^\prime).
$$
Using (1) we get the stated equivalence of functors when evaluated on the subcategory $\CM_{T,\Lambda}^{\varle \lambda+n\chi}$ for any $n$. Now any module that admits a Verma flag lies in one of those categories, hence we proved the fourth claim.  
\end{proof}

For later use we state the following
\begin{lemma}\label{center and antidominant endos} Suppose $T=R$.
The action $Z_{R,\Lambda}\to\End(P^\infty_R(\lambda))$ is an isomorphism.
\end{lemma}
\begin{proof} We have
$$
Z_{R,\Lambda}  =  \bigcap Z_{R_\fp,\Lambda}\subset Z_{R_{(0)},\Lambda},
$$
where the intersection is taken over all prime ideals of height one. Since $P^\infty_R(\lambda)$ is free over $R$, we analogously have
$$
\End(P^\infty_R(\lambda)) = \bigcap \End(P^\infty_R(\lambda)\otimes_R R_\fp)\subset \End(P^\infty_R(\lambda)\otimes_R R_{(0)}).
$$
By Remark \ref{remarks on multiplicity result}, (\ref{remarks on multiplicity result - splitting}) the localizations of $P^\infty_R(\lambda)$ split into the direct sum of the resp. antidominant projectives. Moreover the evaluation maps $Z\to\End(P)$ are compatible with localizations. Since, after localization, we deal with generic or subgeneric situations, the lemma follows from Lemma \ref{calculation of the center - subgeneric case}.
\end{proof}

Now we want to extend the definition of  $\DV$ to the positive level. Let  $\Lambda\in\fhd/_{\textstyle{\sim_T}}$ be an equivalence class outside the critical hyperplanes  and of positive level and let $\lambda\in\Lambda$ be the dominant element. We use the tilting equivalence $t\colon \CM_{T,\Lambda}\stackrel{\sim}\to \CM^{opp}_{T,t(\Lambda)}$. Note that $-2\rho-\lambda$ is an {\em antidominant} weight in $t(\Lambda)$, hence we can define 
\begin{eqnarray*}
\tilde\DV:=\DV_{T,t(\Lambda)}\circ t\colon\CM_{T,\Lambda} & \to & (Z_{T,t(\Lambda)}\catmod)^{opp}\\
M & \mapsto & \Hom_{\CO}(P_T^\infty(-\lambda-2\rho), t(M)).
\end{eqnarray*}
By Lemma \ref{tilted center} the tilting equivalence induces an isomorphism $Z_{T,t(\Lambda)}\cong Z_{T,\Lambda}$, so we can consider $\tilde\DV$ as a functor from $\CM_{T,\Lambda}$ to $(Z_{T,\Lambda}\catmod)^{opp}$. Moreover, by Proposition \ref{properties of V}, (\ref{properties of V - freeness over T}), its image is contained in  the subcategory of $Z_{T,\Lambda}$--modules which are free of finite rank over $T$. Hence composition with the duality $\star=\Hom_T(\cdot,T)$ provides a functor
\begin{eqnarray*}
\DV=\DV_{T,\Lambda}:=\star\circ\DV_{T,t(\Lambda)}\circ t\colon\CM_{T,\Lambda} & \to & Z_{T,\Lambda}\catmod\\
M & \mapsto & \Hom(P_T^\infty(-2\rho-\lambda), t(M))^\star
\end{eqnarray*}
which transforms short exact sequences to short exact sequences, and the statements (2)--(4) of Proposition \ref{properties of V} carry over. 

\begin{remark} 
If the equivalence class $\Lambda$ happens to have both a dominant and an antidominant weight, then the two definitions of $\DV$ agree up to a non--unique isomorphism which can be defined as follows. 

Let $M\in\CM_{T,\Lambda}$ and let $\lambda$ and $\lambda^\prime$ be the antidominant weights of $\Lambda$ and $t(\Lambda)$. Let $\DV_1 M = \Hom(P_T(\lambda),M)$ and $\DV_2 M = \Hom(P_T(\lambda^\prime), t(M))^\star$. Note that $S_{2\rho}\otimes_U P_T(\lambda^\prime)$ is the dual of a non--split module with a finite Verma flag where every Verma module occurs with multiplicity one. Hence it is the dual of the antidominant projective, which is selfdual. Hence  $S\otimes_U P_T(\lambda^\prime)$ is the antidominant projective. Choose generators $v\in P_{T}(\lambda)$ and $v^\prime\in S\otimes_U P_T(\lambda^\prime)$. We get  maps
\begin{eqnarray*}
\DV_1 M \times(\DV_2 M)^\star &\cong& \Hom(P_T(\lambda),M)\times\Hom(M, t(P_T(\lambda^\prime)))\\
 & \to & \Hom(P_T(\lambda),t(P_T(\lambda^\prime))) \\
& = & \Hom(P_T(\lambda), (S\otimes_U P_T(\lambda^\prime))^\circledast) \\
 & \to & T,
\end{eqnarray*}
where the last map is evalution of the image of $v$ on $v^\prime$. This defines a non--degenerate pairing which is functorial in $M$ and hence provides an isomorphism $\DV_1\cong \DV_2$.
\end{remark}

\subsection{The case of general (non--critical) $\Lambda$}
We return to the assumption that $T=R$. Let $\Lambda,\Lambda^\prime\in\fhd/_{\textstyle{\sim}}$ as in Section \ref{subsection - translation functors} and let $\vartheta_{on}$ and $\vartheta_{out}$ be the translation functors. Suppose $\Lambda$ (and hence also $\Lambda^\prime$) is of negative level and let $\lambda\in\Lambda$ and $\lambda^\prime\in\Lambda^\prime$ be the smallest elements. The isomorphism $P_R^\infty(\lambda)\cong\varprojlim \vartheta_{out}P_R^{\varle\lambda^\prime+n\chi}(\lambda^\prime)$ of Remark \ref{remarks on multiplicity result}, (\ref{remarks on multiplicity result - transferred antidominant projective}), induces, via Lemma \ref{center and antidominant endos}, a map $Z_{R,\Lambda^\prime}\to Z_{R,\Lambda}$ which is canonical.  It can be described as follows. An element $z=(z_w)\in Z_{R,\Lambda^\prime}$ will act as multiplication with $z_w$ on $M_R(w.\lambda^\prime)$ for each $w\in\CW(\Lambda)/\Stab(\lambda^\prime)$. Its image will act as multiplication with $z_w$ on $M_R(wx.\lambda)$ for each $x\in\Stab(\lambda^\prime)$. We deduce that $Z_{R,\Lambda^\prime}\to Z_{R,\Lambda}$ is injective, and in the case of regular $\Lambda$ the image is $Z_{R,\Lambda}^{\Stab(\lambda^\prime)}$, i.e.÷ the $\Stab(\lambda^\prime)$--invariant elements under the operation of $\CW(\Lambda)$ defined in Section \ref{subsection - the center}. Using Lemma \ref{tilted center} we also get a map $Z_{R,\Lambda^\prime}\to Z_{R,\Lambda}$ in the case that $\Lambda$ and $\Lambda^\prime$ are of positive level. Let $\Ind= \Ind_{Z_{R,\Lambda^\prime}}^{Z_{R,\Lambda}}=Z_{R,\Lambda}\otimes_{Z_{R,\Lambda^\prime}}\cdot$ be the induction functor and $\Res=\Res^{Z_{R,\Lambda^\prime}}_{Z_{R,\Lambda}}$ its right adjoint.

The following theorem generalizes the subgeneric situation in Theorem \ref{comb of translation - subgeneric}.
\begin{theorem} \label{comb of translation} Suppose $T=R$ and choose $\Lambda$ and $\Lambda^\prime$ as in Section \ref{subsection - translation functors}.
There are isomorphisms of functors
$$
\DV\circ\vartheta_{on}\cong\Res\circ\DV\colon \CM_{R,\Lambda}\to Z_{R,\Lambda^\prime}\catmod
$$
and
$$
\DV\circ\vartheta_{out}\cong\Ind\circ\DV\colon \CM_{R,\Lambda^\prime}\to Z_{R,\Lambda}\catmod.
$$
\end{theorem}

\begin{proof} We will prove the second statement, the proof of the first one is similar.
 In the following all objects will be free over $R$ and hence they are the intersections of their localizations at prime ideals of height one. Moreover, $\DV$ and the translation functors commute with base change by Theorem \ref{properties of translation functors} and Proposition \ref{properties of V}.  For any $R$--module $N$ and any prime $\fp\subset R$ let $N_\fp$ be the localization of $N$ at $\fp$. Let $M\in\CM_{R,\Lambda^\prime}$. We then have
\begin{eqnarray*}
\DV\circ\vartheta_{on} M & = & \bigcap (\DV\circ\vartheta_{on}M)_\fp \\
& = & \bigcap \DV_{\fp}\circ\vartheta_{\fp,on} M_\fp. \\
\end{eqnarray*}
Note that for any $R$--linear functor $F$ we write $F_\fp$ for its version on the localized category. 
The isomorphisms in Theorem \ref{comb of translation - subgeneric} are compatible with base change and with taking the intersection, hence we get
$$
\bigcap \DV_{\fp}\circ\vartheta_{\fp,on}M_\fp=\bigcap \Ind_\fp\circ\DV_\fp M_\fp=\bigcap Z_{R_\fp,\Lambda}\otimes_{Z_{R_\fp,\Lambda^\prime}}\DV_\fp M_\fp.
$$
Since $M$ has a finite Verma flag, the action of $Z_{R_\fp,\Lambda}$ on $Z_{R_\fp,\Lambda}\otimes_{Z_{R_\fp,\Lambda^\prime}}\DV_\fp M_\fp$ factors over a map $Z_{R_\fp,\Lambda}\to Z^I_{R_\fp,\Lambda}$, where $I\subset\CW(\Lambda)/\Stab(\lambda)$ is a finite subset and $Z^I_{R_\fp,\Lambda}$ is the image of $Z_{R_\fp,\Lambda}\subset\prod_{w\in\CW(\Lambda)/\Stab(\lambda)} R_\fp$ under the projection $\prod_{w\in\CW(\Lambda)/\Stab(\lambda)}R_\fp\to \prod_{w\in I} R_\fp$. For finite $I$ we have $Z_{R,\Lambda}^I\otimes_R R_\fp\cong Z_{R_\fp,\Lambda}^I$. Hence 
\begin{eqnarray*}
\bigcap Z_{R_\fp,\Lambda}\otimes_{Z_{R_\fp,\Lambda^\prime}}\DV_\fp M_\fp & = &  \bigcap Z^I_{R_\fp,\Lambda}\otimes_{Z^I_{R_\fp,\Lambda^\prime}}\DV_\fp M_\fp \\
& = & \bigcap (Z^I_{R,\Lambda}\otimes_{Z^I_{R,\Lambda^\prime}}\DV M)_\fp \\
& = & \bigcap (Z_{R,\Lambda}\otimes_{Z_{R,\Lambda^\prime}}\DV M)_\fp \\
& = & \Ind\circ\DV M,
\end{eqnarray*}
as was to be shown. 
\end{proof} 

The next theorem is the analog of Soergel's Struktursatz. Again we do not assume that $M$ or $M^\prime$ is projective, hence it does not hold in this generality for $T=\DC$. 

\begin{theorem}\label{struktursatz} Let $T=R$ and let $\Lambda\in\fhd/_{\textstyle{\sim}}$ be an equivalence class outside the critical hyperplanes. Then  $\DV=\DV_{R,\Lambda}\colon\CM_{R,\Lambda}\to Z_{R,\Lambda}\catmod$ is fully faithful, i.e.~ 
$$
\Hom(M,M^\prime)\stackrel{\sim}\to\Hom(\DV M,\DV M^\prime)
$$
is an isomorphism for any $M,M^\prime\in\CM_{R,\Lambda}$. 
\end{theorem}

\begin{proof} In the following all intersections are taken over the set of prime ideals $\fp\subset R$ of height one. 
Since $M=\bigcap M_\fp\subset M_{(0)}$ we have
$$
\Hom_{\fg\otimes_\DC R}(M,M^\prime)=\bigcap \Hom_ {\fg\otimes_\DC R_\fp}(M_\fp,M^\prime_\fp)\subset \Hom_ {\fg\otimes_\DC R_{(0)}}(M_{(0)}, M^\prime_{(0)}).
$$
Analogously, by Proposition \ref{properties of V}, $\DV N$ is free over $R$ for any $N\in\CM_{R,\Lambda}$. Since $N$ has a finite Verma flag, the module $(\DV N)_\fp$ carries a natural $Z_{R_\fp,\Lambda}$--action. Hence 
\begin{eqnarray*}
\Hom_{Z_{R,\Lambda}}(\DV M,\DV M^\prime) & = & \bigcap\Hom_{Z_{R_\fp,\Lambda}}((\DV M)_\fp,(\DV M^\prime)_\fp) \\
 & \subset &\Hom_{Z_{R_{(0)},\Lambda}}((\DV M)_{(0)},(\DV M^\prime)_{(0)}).
\end{eqnarray*}
All maps are compatible with localizations. Since the localized situations deal with generic and subgeneric cases we can deduce the theorem from Theorem \ref{struktursatz - subgeneric case}.   
\end{proof}

\section{The combinatorics of category $\CO$}\label{section - combinatorics of O}

Let $\Lambda\in\fhd/_{\textstyle{\sim}}$ be an equivalence class outside the critical hyperplanes. We want to give a description of the categorical structure of $\CO_\Lambda$. Suppose $\Lambda$ is in positive level, i.e.÷ $\Lambda$ contains a dominant weight $\lambda$.  For any $w\in\CW(\Lambda)$  there is a projective cover $P(w.\lambda)$ of $L(w.\lambda)$ in $\CO_\Lambda$ and $\CP_\Lambda:=\{P(w.\lambda)\}_{w\in\CW(\Lambda)}$ is a faithful set of small projectives in the sense of \cite{Mit}. We view $\CP_\Lambda$ as a full subcategory of $\CO_\Lambda$.  Let $\DC\catmod^{\CP^{opp}_\Lambda}$ be the category of additive functors $\CP^{opp}_\Lambda\to\DC\catmod$. Then by a theorem of Freyd (cp.ø\cite{Mit}, Theorem 3.1) 
\begin{eqnarray*}
\CO_\Lambda & \to & \DC\catmod^{\CP^{opp}_\Lambda}\\
M & \mapsto & \Hom(\cdot, M)
\end{eqnarray*}
is an equivalence of $\DC$--categories. Hence, in order to describe the category $\CO_\Lambda$, it is sufficient to describe $\CP_\Lambda$. This amounts to a description of $\Hom(P(w.\lambda), P(w^\prime.\lambda))$ for any pair $w,w^\prime\in\CW(\Lambda)$ together with the composition data. Now $P(w.\lambda)\cong P_R(w.\lambda)\otimes_R \DC$ and by Theorem \ref{structure of projectives} we have to describe $\Hom(P_R(w.\lambda), P_R(w^\prime.\lambda))$ and the composition data. By Theorem \ref{struktursatz} we have to describe the $Z_{T,\Lambda}$--modules $\DV P_R(w.\lambda)$.

Suppose $\Lambda$ is regular and of positive level. Let $\lambda\in\Lambda$ be the dominant element. Then $P_R(\lambda)\cong M_R(\lambda)$ and $\DV P_R(\lambda)\cong Z_{R,\Lambda}/\fm_e$, where $\fm_e$ is the annihilator of $M_R(\lambda)$.  Let $w\in\CW(\Lambda)$  and choose a reduced expression $w=s_{1}\cdots s_{n}$  in the Coxeter system $(\CW(\Lambda),\CS(\Lambda))$. We construct $\DV P_R(w.\lambda)$ inductively on the length of $w$. Let $\Theta_{s}$ be the translation functor corresponding to a reflection $s\in\CS(\Lambda)$. By Theorem \ref{properties of translation functors} the module  $\Theta_{s_n}\cdots \Theta_{s_1} M_R(\lambda)$ is projective and the Verma module $M_R(w.\lambda)$ occurs with multiplicity one in a Verma flag since the expression for $w$ was reduced. Hence $P_R(w.\lambda)$ is the unique direct summand that is not isomorphic to $P_R(w^\prime.\lambda)$ for any $w^\prime$ of smaller length. Then  $\DV P_R(w.\lambda)$ is the unique direct summand of $Z_{R,\Lambda}\otimes_{Z_{R,\Lambda}^{s_{n}}}Z_{R,\Lambda}\cdots\otimes_{Z_{R,\Lambda}^{s_{1}}}Z_{R,\Lambda}/\fm_e$ which is not isomorphic to $\DV P_R(w^\prime.\lambda)$ for  any $w^\prime$ of smaller length than $w$. Hence we get an inductive description of $\DV P_R(w.\lambda)$ inside the category of $Z_{R,\Lambda}$--modules. 

If $\Lambda$ is of positive level and not regular, but $\Stab(\lambda)$ is finite, we can find a regular $\Lambda^\prime$ together with the translation on the wall $\vartheta_{on}\colon \CM_{R,\Lambda^\prime}\to\CM_{R,\Lambda}$. There is an induced map $Z_{R,\Lambda^\prime}\to Z_{R,\Lambda}$ which identifies $Z_{R,\Lambda^\prime}$ with the $\Stab(\lambda)$--invariant elements in $Z_{R,\Lambda}$.
Let $w\in\CW(\Lambda)$ and let $\bar w\in\CW(\Lambda)/\Stab(\lambda)$ be  image. Then $\vartheta_{on} P_R(w.\lambda^\prime)$ splits into $\#\Stab(\lambda)$ copies of $P_R(\bar w.\lambda)$ and, possibly, additional terms $P_R(\bar w^\prime.\lambda)$ with $\bar w^\prime.\lambda>\bar w.\lambda$. Analogously, $\DV \vartheta_{on} P(w.\lambda^\prime)=\Res \DV P(w.\lambda^\prime)$ splits into $\#\Stab(\lambda)$ copies of $\DV P(\bar w.\lambda)$ and, possibly, additional terms $\DV P_R(\bar w^\prime.\lambda)$ with $\bar w^\prime.\lambda>\bar w.\lambda$. Again we get an inductive description of $\DV P_R(\bar w.\lambda)$ for any $\bar w\in\CW(\Lambda)/\Stab(\lambda)$. 
So we arrive at a description of $\CP_\Lambda$, hence of $\CO_\Lambda$, in the case that $\Lambda$ is of positive level. 

\begin{theorem}\label{equivalences of categories}
Let $\fh\subset\fb\subset\fg$ and $\fh^\prime\subset\fb^\prime\subset\fg^\prime$ be two symmetrizable Kac-Moody algebras together with Cartan and Borel subalgebras. Choose   two equivalence classes $\Lambda\in\fhd/_{\textstyle{\sim}}$ and $\Lambda^\prime\in(\fh^\prime)^\star/_{\textstyle{\sim^\prime}}$ outside the critical hyperplanes and let $\CO_\Lambda$ and $\CO^\prime_{\Lambda^\prime}$ be the corresponding blocks. Suppose the following.
\begin{enumerate}
\item There exist $\lambda\in\Lambda$ and $\lambda^\prime\in\Lambda^\prime$ which are either both dominant or both antidominant.
\item There is an isomorphism $(\CW(\Lambda),\CS(\Lambda))\cong(\CW^\prime(\Lambda^\prime),\CS^\prime(\Lambda^\prime))$ of the corresponding Coxeter systems. 
\item This isomorphism induces a bijection $\Stab(\lambda)\cong\Stab(\lambda^\prime)$ an both sets are finite. 
\end{enumerate}
Then there is an equivalence of categories
$$
\CO_\Lambda\cong\CO^\prime_{\Lambda^\prime}.
$$
\end{theorem}
\begin{proof}
The description of $\CP_\Lambda$ above proves the theorem in the case of $\Lambda$ and $\Lambda^\prime$ in positive level. If $\Lambda$ and $\Lambda^\prime$ are in negative level, then $t(\Lambda)$ and $t(\Lambda^\prime)$ are in positive level. Since Verma modules have a categorical characterization as projective covers in truncated categories, any equivalence $\CO_{t(\Lambda)}\cong \CO^\prime_{t(\Lambda^\prime)}$ induces an equivalence $\CM_{t(\Lambda)}\cong \CM^\prime_{t(\Lambda^\prime)}$. Using the tilting functor we get an equivalence $\CM_{\Lambda}\cong\CM^\prime_{\Lambda^\prime}$. It induces an equivalence $\CM_{\Lambda}^{\varle w.\lambda}\cong(\CM^\prime_{\Lambda^\prime})^{\varle w.\lambda^\prime}$ for any $w\in\CW(\Lambda)\cong\CW^\prime(\Lambda^\prime)$. This equivalence must map the indecomposable projective objects in $\CO_\Lambda^{\varle w.\lambda}$ to the indecomposable projective objects of $(\CO_{\Lambda^\prime}^\prime)^{\varle w.\lambda^\prime}$, hence it induces an equivalence $\CO_\Lambda^{\varle w.\lambda}\cong (\CO^\prime_{\Lambda^\prime})^{\varle w.\lambda^\prime}$. Since a $\fg$-module is in $\CO$ if and only if any finitely generated submodule is in some truncated subcategory, we get an equivalence $\CO_\Lambda\cong \CO^\prime_{\Lambda^\prime}$.
\end{proof}

\section{An application}
Let $\Lambda\in\fhd_{\textstyle{\sim}}$ be a regular equivalence class outside the critical hyperplanes and suppose that the Coxeter system $(\CW(\Lambda),\CS(\Lambda))$ is of finite or affine type. According to Theorem \ref{equivalences of categories} the block $\CO_\Lambda$ is equivalent to a regular integral block $\CO^\prime_{\Lambda^\prime}$ of the finite of affine Kac-Moody algebra $\fg^\prime$ for the Coxeter system $(\CW^\prime,\CS^\prime)\cong(\CW(\Lambda),\CS(\Lambda))$ (note that regular integral equivalence classes exist in those cases). Let $\lambda\in\Lambda$ and $\lambda^\prime\in\Lambda^\prime$ be both dominant or both antidominant. Under the above  equivalence  the class of the Verma module $M(w.\lambda)$ must map to the class of $M(w.\lambda^\prime)$ and the class of the simple module $L(w.\lambda)$ must map to the class of $L(w.\lambda^\prime)$.  Hence $[M(x.\lambda): L(y.\lambda)]=[M(x.\lambda^\prime):L(y.\lambda^\prime)]$. So 
we reduced the following statement to the integral case, where it was proved  by Beilinson-Bernstein \cite{BB} and Brylinski-Kashiwara \cite{BK} in the finite situation and by Kashiwara \cite{Kas} and Kashiwara-Tanisaki \cite{KT1} in the dominant affine and Kashiwara-Tanisaki \cite{KT2} in the antidominant affine case.
\begin{theorem}\label{Kazhdan-Lusztig conjecture} Let $\fg$ be a symmetrizable Kac-Moody algebra and $\Lambda\in\fhd/_{\textstyle{\sim}}$ an equivalence class outside the critical hyperplanes and let $\lambda\in\Lambda$ be dominant or antidominant. Suppose that $(\CW(\Lambda),\CS(\Lambda))$ is of finite or affine type and  that  $\Lambda$ is regular, i.e.ø $\Stab(\lambda)=\{e\}$. 
 Then the Kazhdan-Lusztig conjecture holds, i.e.ø
$$
\ch\, L(w.\lambda) =
 \sum_{y\ge w}(-1)^{l(y)-l(w)}Q_{w,y}(1)\,\ch\, M(y.\lambda)
$$ 
if $\lambda$ is dominant and, if $\lambda$ is antidominant,
$$
\ch\, L(w.\lambda) = \sum_{y\le w}(-1)^{l(w)-l(y)}P_{y,w}(1)\,\ch\, M(y.\lambda),
$$
where $P_{y,w}$ and $Q_{w,y}\in\DZ[v,v^{-1}]$ are the Kazhdan-Lusztig polynomial and the inverse Kazhdan-Lusztig polynomial for the Coxeter system $(\CW(\Lambda),\CS(\Lambda))$.
\end{theorem}

 In fact, Kashiwara \cite{Kas} and Kashiwara-Tanisaki \cite{KT1} proved the conjecture for any symmetrizable Kac-Moody algebra in the integral dominant case. Then Kashiwara-Tanisaki proved it subsequently for the integral antidominant case and affine algebras \cite{KT2}, for the rational antidominant case and affine algebras \cite{KT3}, in the  rational dominant case and arbitrary symmetrizable Kac-Moody algebras \cite{KT4}, and for the arbitrary (non--critical) dominant or antidominant case and affine algebras \cite{KT5}.

\end{document}